\documentclass[12pt,a4paper]{article}
\usepackage{amssymb,amsmath}
\usepackage{amsfonts}
\usepackage{centernot}
\usepackage{mathtools}
\usepackage{ stmaryrd }
\usepackage{a4}
\newcommand{\para}{\par\vspace{.25cm}}

\newtheorem{theorem}{Theorem}
\newtheorem{lemma}{Lemma}
\newtheorem{cor}{Corollary}
\newtheorem{remark}{Remark}
\usepackage{multirow}

\begin{document}
\baselineskip 18pt

\title{\bf On the index of a free abelian subgroup in the group of central units of an integral group ring}
\author{ Gurmeet K. Bakshi and Sugandha Maheshwary{\footnote {Research supported by CSIR, India, is gratefully acknowledged} \footnote{Corresponding author}} \\ {\em \small Centre for Advanced Study in
Mathematics,}\\
{\em \small Panjab University, Chandigarh 160014, India}\\{\em
\small email: gkbakshi@pu.ac.in and msugandha.87@gmail.com } }
\date{}
{\maketitle}

\begin{abstract}\noindent
{Let $\mathcal{Z}(\mathcal{U}(\mathbb{Z}[G]))$ denote the group of central units in the integral group ring $\mathbb{Z}[G]$ of a finite group $G$. A bound on the index of the subgroup generated by  a virtual basis in $\mathcal{Z}(\mathcal{U}(\mathbb{Z}[G]))$ is computed for a class of strongly monomial groups. The result is illustrated with application to the groups of order $p^{n}$, $p$ prime, $n \leq 4$. The rank of $\mathcal{Z}(\mathcal{U}(\mathbb{Z}[G]))$ and the Wedderburn decomposition of the rational group algebra of these $p$-groups have also been obtained.}
\end{abstract}\vspace{.25cm}
 \vspace{.25cm} 

\section{Introduction} Let $\mathcal{U}(\mathbb{Z}[G])$ denote the unit group of the integral group ring $\mathbb{Z}[G]$ of a finite group $G$. The centre of $\mathcal{U}(\mathbb{Z}[G])$ is denoted by $\mathcal{Z}(\mathcal{U}(\mathbb{Z}[G]))$. It is well known that $\mathcal{Z}(\mathcal{U}(\mathbb{Z}[G]))= \pm \mathcal{Z}(G) \times A$, where $A$ is a free abelian subgroup of $\mathcal{Z}(\mathcal{U}(\mathbb{Z}[G]))$ of finite rank. In order to study $\mathcal{Z}(\mathcal{U}(\mathbb{Z}[G]))$, a multiplicatively independent subset of such a subgroup $A$, i.e., a $\mathbb{Z}$-basis for such a free $\mathbb{Z}$-module $A$, is of importance, and is known only for a few groups (\cite{aleev,aleev2,Ferr,Li}, see also \cite{milies}, Examples 8.3.11 and 8.3.12). However, other papers deal with determining a  virtual basis of  $\mathcal{Z}(\mathcal{U}(\mathbb{Z}[G]))$, i.e., a multiplicatively independent subset of $\mathcal{Z}(\mathcal{U}(\mathbb{Z}[G]))$ which generates a subgroup of finite index in $\mathcal{Z}(\mathcal{U}(\mathbb{Z}[G]))$ (see e.g. \cite{Ferraz,hoech,hoech2,hoech3,jesp1,jesp3,jesp4,jesp,jesp2,parm}).
 \par Analogous to well known cyclotomic units in cyclotomic fields, Bass \cite{Bass} constructed units, so called \emph{Bass cyclic units}, which generate  a subgroup of finite index in $\mathcal{U}(\mathbb{Z}[G])$, when $G$ is cyclic. A virtual basis consisting of certain Bass cyclic units was also given by Bass. Generalizing the notion of Bass cyclic units, Jespers et al \cite{jesp3} defined \emph{generalized Bass units} and have shown that the group generated by these units contains a subgroup of finite index in $\mathcal{Z}(\mathcal{U}(\mathbb{Z}[G]))$ for an arbitrary strongly monomial group $G$. Recently, for a class of groups properly contained in finite strongly monomial groups, Jespers et al \cite{jesp} provided a subset, denoted by $\mathcal{B}(G)$ (say), of the group generated by generalized Bass units, which forms a virtual basis of $\mathcal{Z}(\mathcal{U}(\mathbb{Z}[G]))$. \par In this paper, we determine a bound on the index of the subgroup generated by $\mathcal{B}(G)$  in $\mathcal{Z}(\mathcal{U}(\mathbb{Z}[G]))$ for the same class of groups as considered in \cite{jesp} (Theorem \ref{t1}). Our result is based on the ideas contained in \cite{jesp} and Kummer's work (see \cite{kumm}, Theorem 8.2) on the index of cyclotomic units. Further in \cite{jesp}, Jespers et al have provided the rank of $\mathcal{Z}(\mathcal{U}(\mathbb{Z}[G]))$ in terms of strong Shoda pairs of $G$, when $G$ is a strongly monomial group. In section 4, we compute a complete and irredundant set of strong Shoda pairs of the non abelian groups of order $p^{n}$, $p$ prime, $n \leq 4$, and provide, in terms of $p$, the rank of $\mathcal{Z}(\mathcal{U}(\mathbb{Z}[G]))$ of these $p$-groups  along with the Wedderburn decomposition of their rational group algebras. We also illustrate Theorem \ref{t1} for some of these groups.

\section{ Notation and Preliminaries}
We begin by fixing some notation.\\

$ \begin{array}{lll}
    G  &{\rm a ~finite~ group} \vspace{.1cm}\\
    |g| &{\rm the~ order ~of~the ~ element}~ g~{\rm in}~ G \vspace{.1cm} \\
    g^t & t^{-1}gt,~ g, t \in G\vspace{.1cm} \\
\langle X\rangle & {\rm the ~subgroup~ generated ~by~ the ~subset}~ X ~{\rm of}~ G\vspace{.1cm}\\
|X| & {\rm the~ cardinality ~of ~the ~set}~ X\vspace{.1cm}\\
 K \leq G  & K {\rm~is~ a ~subgroup ~of} ~G\vspace{.1cm}\\
    K \unlhd G &  K {\rm ~is~ a~ normal~ subgroup~ of}~ G\vspace{.1cm}\\
    \left[G:K\right] &  {\rm the~ index~ of~ the~ subgroup}~ K ~{\rm in}~  G\vspace{.1cm}\\
      N_{G}(K) & {\rm the~ normalizer~ of}~ K ~{\rm in}~ G\vspace{.1cm}\\
\operatorname{core}(K) & \bigcap_{x\in G}xKx^{-1} ,~{\rm the~ largest~ normal~ subgroup~ of ~}G ~{\rm contained~in}~ K\vspace{.1cm}\\\end{array}
$

$\begin{array}{lll}
\hat{K}  &\frac{1}{|K|} { \sum_{k \in K}k}\vspace{.1cm}\\
\mathcal{M}(G/K) & {\rm the ~set ~of~ minimal~ normal~ subgroups~ of}~ G~ {\rm containing}~ K {\rm ~properly} \vspace{.1cm} \\
\varepsilon (H,K) & \begin{cases} \hat{H}, & {\rm if}~ H = K; \\ \prod_{M/K \in \mathcal{M}(H/K)} (\hat{K}-\hat{M}), &  {\rm otherwise,~where~}K\unlhd H \leq G \end{cases}\vspace{.1cm}\\
e(G,H,K) & {\rm the~ sum ~of ~all ~the~ distinct}~ G{\rm-conjugates~ of}~ \varepsilon (H, K)\vspace{.1cm}\\
  \varphi & {\rm Euler's~ phi~ function}\vspace{.1cm}\\
  \mathbb{Z}/n\mathbb{Z} & {\rm the ~ring ~of ~integers ~modulo}~ n,~ n \geq 1 \vspace{.1cm}\\
  \zeta_{n}& {\rm a ~primitive}~ n^{th} {\rm~ root~ of~ unity ~in ~the~ field~ of~ complex~ numbers}\vspace{.1cm}\\
  \operatorname{Gal}(\mathbb{Q}(\zeta_{n})/\mathbb{Q})&{\rm the~Galois ~group ~of~ the~ cyclotomic~ field }~ \mathbb{Q}(\zeta_{n})~ {\rm over}~ \mathbb{Q}\vspace{.1cm}\\
  h_{n}^{+} & {\rm the~ class~ number~ of ~the ~maximal ~real~ subfield ~of}~ \mathbb{Q}(\zeta_{n})\vspace{.1cm}\\
l.c.m.(k,n)& {\rm the~ least ~common~ multiple~ of~the~integers} ~k~ {\rm and}~ n\vspace{.1cm}\\
(k,n)& {\rm the~ greatest ~common~ divisor~ of~the~integers} ~k~ {\rm and}~ n\vspace{.1cm}\\
  o_{n}(k)& {\rm the~ multiplicative~ order~ of}~ k~{\rm modulo}~ n, {\rm~where~} (k,n)=1\vspace{.1cm}\\
\eta_{k}(\zeta_{n})&\begin{cases}  1,& {\rm if} ~n=1; \\ 1+\zeta_{n}+ \zeta_{n}^{2}+...+\zeta_{n}^{k-1}, &{\rm if}~ n>1,~{\rm where~} k\geq 1\end{cases} \vspace{.1cm}\\
\mathcal{U}(R)& {\rm the ~unit~ group~ of ~ the~ ring}~ R \vspace{.1cm}\\
M_{n}(R)& {\rm the~ ring ~of}~ n \times n~ {\rm matrices~ over~ the~ ring}~R,~n \geq 1\vspace{.1cm}\\
M_{n}(R)^{(s)}& M_{n}(R)\oplus M_{n}(R)\oplus...\oplus M_{n}(R),~{\rm direct ~sum ~of}~ s~{\rm copies},~ s\geq 1\vspace{.1cm}\\
I_{n} & {\rm the}~n \times n~ {\rm identity ~matrix}\vspace{.1cm}\\

\end{array}$

\vspace{1cm}

\par A {\it strong Shoda pair} (\cite{Oli}, Definition 3.1) of $G$ is a pair $(H,
K)$ of subgroups of $G$ with the properties that \begin{quote}  (i) $K \unlhd H \unlhd N_{G}(K)$;\\(ii) $H/K$ is cyclic and a maximal abelian subgroup of $N_{G}(K)/K$; \\ (iii) the distinct $G$-conjugates of $\varepsilon(H,K)$ are mutually orthogonal. \end{quote}
Note that $(G,G)$ is always a strong Shoda pair of $G$.\\

If $(H,K)$ is a strong Shoda pair of $G$, then $e(G,H,K)$ is a primitive central idempotent of the rational group algebra $\mathbb{Q}[G]$ (\cite{Oli}, Proposition 3.3). A group $G$ is called \emph{strongly monomial} if every primitive central idempotent of $\mathbb{Q}[G]$ is of the form $e(G,H,K)$ for some strong Shoda pair $(H,K)$ of $G$.\\

Two strong Shoda pairs $(H_{1}, K_{1})$ and $(H_{2}, K_{2})$ of $G$ are said to be {\it equivalent} if $e(G, H_{1}, K_{1})$ = $e(G, H_{2}, K_{2})$. A complete set of representatives of distinct equivalence classes of strong Shoda pairs of $G$ is called {\it a complete irredundant set of strong Shoda  pairs} of $G$. In case $G$ is strongly monomial, one can calculate the primitive central idempotents of $\mathbb{Q}[G]$ from a complete irredundant  set of strong Shoda pairs of $G$.\\

Recall that a group $G$ is called \emph{normally monomial} if every complex irreducible character of $G$ is induced from a linear character of a normal subgroup of $G$. Theorem \ref{t3}, as stated below, provides an algorithm to determine a complete irredundant set of strong Shoda pairs of a normally monomial group $G$ and also, in particular, yields that a normally monomial group is strongly monomial.\\

 \par Let $\mathcal{N}$ be the set of all the distinct normal subgroups of a finite group $G$. For $N \in \mathcal{N}$, set \vspace{.2cm}\\ $\begin{array}{lll} A_{N}: & {\rm a ~normal~ subgroup~ of~} G {\rm~ containing~} N {\rm ~such~ that~ } A_{N}/N {\rm ~is ~an ~abelian}\\ &
  { \rm normal~  subgroup~ of~ maximal~ order~ in~} G/N. \vspace{.2cm}\\
  \mathcal{D}_{N}: & {\rm the ~set ~of ~all ~subgroups~} D ~{\rm  of} ~A_{N}~ {\rm containing ~} N~{\rm  such ~that~}
   \operatorname{core}(D)=N, \\&  A_{N}/D ~{\rm    is ~ cyclic ~  and ~  is ~a ~ maximal~ abelian ~ subgroup~  of~}   N_{G}(D)/D. \vspace{.2cm}\\ \mathcal{T}_{N}: &  {\rm a ~set~ of ~representatives~ of~ } \mathcal{D}_{N} {\rm ~  under ~the~ equivalence~ relation~ defined~ by} \\ & {\rm
    conjugacy~ of~ subgroups~in~} G. \vspace{.2cm} \\  \mathcal{S}_{N}: & \{( A_{N},D)~|~ D \in \mathcal{T}_{N}\}.\end{array}$ \vspace{.2cm}\para     Note that if $N \in \mathcal{N}$ is such that $G/N$ is abelian, then, by (\cite{BM}, Eq.(1)), \begin{equation}\label{s1} \mathcal{S}_{N} = \begin{cases} \{(G, N)\}, & {\rm if} ~G/N~{\rm cyclic;}\\ \emptyset, & {\rm otherwise.}
    \end{cases} \end{equation}

    \begin{theorem}{\rm(\cite{BM}, Theorem 1, Corollaries 1 and 2)} \label{t3} The following statements are equivalent:

     \begin{quote}
 {\rm(i)} $G$ is normally monomial; \\{\rm(ii)} $\mathcal{S}(G):=\bigcup_{N \in \mathcal{N} }\mathcal{S}_{N} $ is a complete irredundant set of strong Shoda pairs of $G$;\\ {\rm(iii)} $\{ e(G, A_{N}, D)\, |\,  (A_{N},D)  \in S_{N},~ N \in \mathcal{N} \}$ is a complete set of primitive central idempotents of $\mathbb{Q}[G]$;\\{\rm(iv)} $ |G| =  \sum_{N\in \mathcal{N}}\sum_{D \in \mathcal{D}_{N}}[G:A_{N}]\varphi([A_{N}:D])$. \\

 \end{quote}
 \end{theorem}
 \vspace{0.3cm}

Let $n \geq 1$ and let $k$ be an integer coprime to $n$. Then, $\eta_{k}(\zeta_{n})$ is a unit of $\mathbb{Z}[\zeta_{n}]$.
 The units of the form $\eta_{k}(\zeta_{n}^{j})$ with integers $j,k$ and $n$ such that $(k,n)=1$ are called \emph{cyclotomic units} of $\mathbb{Q}(\zeta_{n})$.\\

Let $g\in G$ and let $k,~m$ be positive integers such that $k^{m}\equiv 1({\rm mod}~n)$, where $n=|g|$. Then, $$u_{k,m}(g)=(1+g+...+g^{k-1})^{m} + \frac{1-k^{m}}{n}(1+g+...+g^{n-1})$$ is a unit in the integral group ring $\mathbb{Z}[G]$. The units in $\mathbb{Z}[G]$ of this form are called \emph{Bass cyclic units}(see \cite{units}, (10.3)). \\

Next, we recall the definition of \emph{generalized Bass units} of $\mathbb{Z}[G]$ introduced by Jespers et al \cite{jesp3}. For  $M \unlhd G$, $g \in G$, and positive integers $k$, $m$   such that $k^{m}\equiv1 ({\rm mod}~ |g|)$, let $$u_{k,m}(1-\hat{M}+g\hat{M})= 1-\hat{M}+u_{k,m}(g)\hat{M}.$$ This element is a unit in $\mathbb{Z}[G](1-\hat{M})+\mathbb{Z}[G]\hat{M}$. As both $\mathbb{Z}[G](1-\hat{M})+\mathbb{Z}[G]\hat{M}$ and $\mathbb{Z}[G]$ are orders in $\mathbb{Q}[G]$, there is a positive integer $n_{g,M}$ such that
   \begin{equation}\label{eq2}
    (u_{k,m}(1-\hat{M}+g\hat{M}))^{n_{g,M}}\in \mathcal{U}(\mathbb{Z}[G]).
   \end{equation}
 Suppose $n_{G,M}$ is the minimal positive integer satisfying Eq.(\ref{eq2}) for all $g \in G$. Then, the element $$  (u_{k,m}(1-\hat{M}+g\hat{M}))^{n_{G,M}}=1-\hat{M}+u_{k,mn_{G,M}}(g)\hat{M}$$

   \noindent is called the \emph{generalized Bass unit} of $\mathbb{Z}[G]$ based on $g$ and $M$ with parameters $k$ and $m$. Observe that $n_{G,M}=1$, if $M$ is trivial i.e., $M=\langle1\rangle$ or $G$.

\begin{remark}\label{r1}For a non trivial normal subgroup $M$ of $G$, using Lemma 3.1 of {\rm\cite{Gonc}}, it may be noted that $n_{G,M}=1$, if every $g\in G\setminus M$ is of order 2; otherwise, $n_{G,M}$ is the minimal positive integer satisfying Eq.{\rm(\ref{eq2})} for all $g\in G\setminus M$ with $|g|>2$ and integers $k,~m$ such that $1<k<|g|$, $(k,|g|)=1$ and $m=o_{|g|}(k).$\end{remark}

\par Let $G$ be a strongly monomial group such that there is a complete irredundant set $\{(H_{i},K_{i})~|~ 1\leq i\leq m\}$ of strong Shoda pairs of $G$ with the property that each $[H_{i}:K_{i}]$ is a prime power, say $p_{i}^{n_{i}}$. Assume that $(H_{1},K_{1})=(G,G)$. For such a group $G$, we recall the virtual basis of $\mathcal{Z}(\mathcal{U}(\mathbb{Z}[G]))$ provided by Jespers et al \cite{jesp}. For $1 \leq i \leq m$, we adopt the following notation:

$$\begin{array}{ll}
\varepsilon_{i}& := \varepsilon(H_{i},K_{i})\\
e_{i} &:= e(G, H_{i}, K_{i})\\
  \left[H_{i}:K_{i}\right]&:=p_{i}^{n_{i}}, p_{i} ~{\rm prime},~ n_{i} \geq 0~(n_{i}=0 ~{\rm only~if}~i=1)\\
  g_{i}K_{i} & := {\rm a~ generator ~of~ the~ cyclic~ group} ~H_{i}/K_{i} \\
  L_{j}^{(i)}&:=\langle g_{i}^{p_{i}^{n_{i}-j}}, K_{i} \rangle ,~ 0\leq j \leq n_{i}\\
  N_{i}&:=N_{G}(K_{i})\\
 m_{i}&:= [G: N_{i}]\\
 T_{i}&:= {\rm a ~right~ transversal~ of}~ N_{i}~ {\rm in}~ G
 \end{array}$$\\

  \noindent For $2 \leq i \leq m$, consider the action of $N_{i}/H_{i}$ on $\mathbb{Q}(\zeta_{p_{i}^{n_{i}}})$ given by the map
\begin{eqnarray}
\nonumber
  N_{i}/H_{i} &\longrightarrow& \operatorname{Gal}(\mathbb{Q}(\zeta_{p_{i}^{n_{i}}})/\mathbb{Q}) \nonumber\\
  n_{i}H_{i} & \longmapsto & \alpha_{n_{i}H_{i}},
\end{eqnarray}

\noindent where $\alpha_{n_{i}H_{i}}(\zeta_{p_{i}^{n_{i}}})=\zeta_{p_{i}^{n_{i}}}^{j}$, if $n_{i}^{-1}g_{i}n_{i}K_{i}=g_{i}^{j}K_{i}.$ As $H_{i}/K_{i}$ is a maximal abelian subgroup of $N_{i}/K_{i}$, it turns out that the above action is faithful. Hence, $ N_{i}/H_{i}$ is isomorphic to a subgroup of $\operatorname{Gal}(\mathbb{Q}(\zeta_{p_{i}^{n_{i}}})/\mathbb{Q})\cong \mathcal{U}(\mathbb{Z}/p_{i}^{n_{i}}\mathbb{Z})$. For the convenience of notation, we regard $N_{i}/H_{i}$ as a subgroup of $\operatorname{Gal}(\mathbb{Q}(\zeta_{p_{i}^{n_{i}}})/\mathbb{Q})$ and that of $\mathcal{U}(\mathbb{Z}/p_{i}^{n_{i}}\mathbb{Z})$. (Notice that $N_{i}/H_{i}$ can be regarded as a subgroup of $\mathcal{U}(\mathbb{Z}/\left[H_{i}:K_{i}\right]\mathbb{Z})$, even if $\left[H_{i}:K_{i}\right]$ is not a prime power.) With this identification, $ N_{i}/H_{i}$ is equal to either $\langle \phi_{r_{i}}\rangle$ or $\langle \phi_{r_{i}}\rangle \times \langle \phi_{-1}\rangle$ (resp. $\langle r_{i}\rangle$ or $\langle r_{i}\rangle \times \langle-1\rangle$) for some $r_{i} \in \mathcal{U}(\mathbb{Z}/p_{i}^{n_{i}}\mathbb{Z}),$ where $\phi_{r_{i}}$ denotes the automorphism of $\mathbb{Q}(\zeta_{p_{i}^{n_{i}}})$ which sends $\zeta_{p_{i}^{n_{i}}}$ to $\zeta_{p_{i}^{n_{i}}}^{r_{i}}.$ The later case arises only if $p_{i}=2$ and $n_{i}\geq 3 $. Set

\begin{equation}\label{eqn 5}
    d_{i}:= \begin{cases}1,~{\rm if} -1 \in \langle r_{i}\rangle; \\
                     2, ~{\rm otherwise}.\end{cases}
\end{equation}

\noindent and

\begin{equation}\label{eqn 7}
    o_{i}:= \left\{
              \begin{array}{ll}
                 4, & \hbox{if $p_{i}=2$, $N_{i}/H_{i}=\langle r_{i}\rangle$, $r_{i}\equiv1({\rm mod}~ 4)$, $n_{i}\geq 2$;} \\
                6, & \hbox{if $p_{i}=3$, $N_{i}/H_{i}=\langle r_{i}\rangle$, $r_{i}\equiv1({\rm mod} ~3)$;} \\
               2, & \hbox{otherwise.}
              \end{array}
            \right.
\end{equation}\\
\noindent Further, choose a subset  $I_{i}$  of $\{k ~|~ 1 \leq k \leq \frac{{p_{i}}^{n_{i}}}{2}, (k,p_{i})=1 \} $ containing 1, which forms a set of representatives of $\mathcal{U}(\mathbb{Z}/p_{i}^{n_{i}}\mathbb{Z})$ modulo $\langle N_{i}/H_{i}, -1\rangle.$ We extend the notation by setting $I_{1}=\{1\}$, in view of the trivial action of the identity group $N_{1}/H_{1}$ on $\mathbb{Q}(\zeta_{1})=\mathbb{Q}$.
\para Let $k$ be a positive integer coprime to $p_{i}$ and let $r$ be an arbitrary integer. For $0\leq j\leq s \leq n_{i}$, consider the following products of generalized Bass units of $\mathbb{Z}[H_{i}]$, defined recursively:
 $$\begin{array}{ll}
    c_{s}^{s}(H_{i},K_{i},k,r) & = 1 \\\end{array}$$

 \noindent and, for $0\leq j \leq s-1$,
  $$\begin{array}{rl}
c_{j}^{s}(H_{i},K_{i},k,r)~ =& \big(\prod_{h\in L_{j}^{(i)}}u_{k,o_{p_{i}^{n_{i}}}(k)n_{H_{i},K_{i}}}(g_{i}^{rp_{i}^{n_{i}-s}}h\hat{K_{i}}+1-\hat{K_{i}})\big)^{p_{i}^{s-j-1}} \times \vspace{.2cm}\\
     & ~~~~~~~\big(\prod_{l=j+1}^{s-1}c_{l}^{s}(H_{i},K_{i},k,r)^{-1}))(\prod_{l=0}^{j-1}c_{l}^{s+l-j}(H_{i},K_{i},k,r)^{-1}\big),
  \end{array}$$
\noindent where the empty products equal 1.\\

\noindent Define

   $$\begin{array}{lll}
        & B(H_{i},K_{i}) & :=\{\prod_{x\in N_{i}/H_{i}}c_{0}^{n_{i}}(H_{i},K_{i},k,x)~|~ k \in I_{i} \setminus \{1\} \}, \vspace{0.2cm}\\
        & \mathcal{B}(H_{i},K_{i}) & := \{ \prod_{t\in T_{i}}u^{t}~ | ~u \in B(H_{i},K_{i})\}, \vspace{0.2cm}\\

     \end{array}
   $$\\
   and
\begin{equation}\label{eq1}
    \mathcal{B}(G):= \bigcup_{i=1}^{m}B(H_{i},K_{i}).
\end{equation}

\vspace{0.4cm}
 \noindent Jespers et al (\cite{jesp}, Theorem 3.5) proved that $ \mathcal{B}(G)$ is a \emph{virtual basis of} $\mathcal{Z}(\mathcal{U}(\mathbb{Z}[G]))$.

 \section{A bound on the index of $\langle \mathcal{B}(G) \rangle$ in $\mathcal{Z}(\mathcal{U}(\mathbb{Z}[G]))$}

\par In this section, we continue with the notation developed in section 2.

\begin{theorem}\label{t1} Let $G$ be a strongly monomial group and let $\{(H_{i},K_{i})~|~ 1\leq i\leq m\}$ be a complete irredundant set  of strong Shoda pairs of $G$ with $(H_{1},K_{1})=(G,G)$. For $2 \leq i \leq m$, let $I_{i}$ be a subset of $\{k ~|~ 1 \leq k \leq  \frac{[H_{i}:K_{i}]}{2}, (k,[H_{i}:K_{i}])=1 \} $ containing 1, which forms a set of representatives of $\mathcal{U}(\mathbb{Z}/[H_{i}:K_{i}]\mathbb{Z})$ modulo $\langle N_{i}/H_{i}, -1\rangle$, where $N_{i}=N_{G}(K_{i})$. Set $I_{1}=\{1\}$.
\begin{description}
  \item[(i)] The rank of $\mathcal{Z}(\mathcal{U}(\mathbb{Z}[G]))=0$ if and only if $|I_{i}|=1$ for all $i$, $1\leq i\leq m$, and in this case, $\mathcal{Z}(\mathcal{U}(\mathbb{Z}[G]))=\pm \mathcal{Z}(G).$
  \item[(ii)] In addition, if $[H_{i}:K_{i}]$ is a prime power, say $p_{i}^{n_{i}}$, for all $i$, $1\leq i\leq m$, and $\mathcal{B}(G)$ is the virtual basis of $\mathcal{Z}(\mathcal{U}(\mathbb{Z}[G]))$ as defined in Eq.{\rm (\ref{eq1})}, then,

$$[\mathcal{Z}(\mathcal{U}(\mathbb{Z}[G])):\langle \mathcal{B}(G)\rangle] \leq 2\prod_{\stackrel{i=2}{|I_{i}|=1}
       }^{m}o_{i}\prod_{\stackrel{i=2}{|I_{i}|\neq 1}}^{m} h_{p_{i}^{n_{i}}}^{+}l_{i}p_{i}^{n_{i}-1}\mathfrak{o}_{i}(l_{i}^{d_{i}-1}[N_{i} :H_{i}])^{|I_{i}|-1},$$

\noindent where  $l_{i}=l.c.m.(2,p_{i})$; $\mathfrak{o}_{i}=\hspace{-0.3cm}{\displaystyle\prod_{\stackrel{1< k < \frac{p_{i}^{n_{i}}}{2}}{(k,p_{i})=1}}o_{p_{i}^{n_{i}}}(k)p_{i}^{n_{i}-1}n_{H_{i},K_{i}}}$; $ d_{i}$ and $o_{i}$ are as defined in Eq.{\rm(\ref{eqn 5})} and Eq.{\rm(\ref{eqn 7})} respectively.
\end{description}

\end{theorem}

\vspace{0.5cm}

\noindent We first prove the following:

\begin{lemma}\label{l1} Let $G$ be as in Theorem \ref{t1}. Let $\mathcal{A}(H_{i},K_{i})=\mathcal{Z}(1-e_{i}+\mathcal{U}(\mathbb{Z}[G]e_{i}))$ and $A(H_{i}, K_{i})=\mathcal{Z}(1-\varepsilon_{i}+\mathcal{U}(\mathbb{Z}[N_{i}]\varepsilon_{i}))$, where $1 \leq i \leq m$. Then,
$$[\mathcal{A}(H_{i},K_{i}):\langle \mathcal{B}(H_{i},K_{i}) \rangle] = [A(H_{i},K_{i}): \langle B(H_{i},K_{i}) \rangle].$$

\end{lemma}

\noindent{\bf Proof.} Let $\{ t_{j} ~ |~ 1 \leq j \leq m_{i}\}$ be a right transversal of $N_{i}$ in $G$ with $t_{1}=1$. For $ \alpha \in  \mathbb{Q}[G]e_{i}$ and integers $r$ and $s$ such that $ 1 \leq r, s\leq m_{i}$, let $\alpha_{rs} = \varepsilon_{i}t_{r}\alpha  t_{s}^{-1} \varepsilon_{i}$. We notice that $\alpha_{rs} \in \mathbb{Q}[N_{i}]\varepsilon_{i}$. To see this, write $\alpha = (\sum_{g \in G}\alpha_{g}g)e_{i}$ with $\alpha_{g} \in \mathbb{Q}$. Then $\alpha_{rs}  = \sum_{g \in G} \alpha_{g}\varepsilon_{i}t_{r}gt_{s}^{-1}\varepsilon_{i}$. By (\cite{Oli}, Proposition 3.3), the centralizer of $\varepsilon_{i}$ in $G$ equals $N_{i}$. Therefore, if $t_{r}gt_{s}^{-1} \not \in N_{i}$, then $\varepsilon_{i}t_{r}gt_{s}^{-1}\varepsilon_{i}= t_{r}gt_{s}^{-1}\varepsilon_{i}^{t_{r}gt_{s}^{-1}}\varepsilon_{i}=0.$ Also, if $t_{r}gt_{s}^{-1} \in N_{i}$, then $\varepsilon_{i}t_{r}gt_{s}^{-1}\varepsilon_{i}=t_{r}gt_{s}^{-1}\varepsilon_{i}.$  Consequently, $\alpha_{rs}  = \sum_{ t_{r}gt_{s}^{-1} \in N_{i}} \alpha_{g} t_{r}gt_{s}^{-1}\varepsilon_{i} \in \mathbb{Q}[N_{i}]\varepsilon_{i}$.\\

 \noindent Now consider the map\vspace{-0.3cm}$$\theta_{i}:\mathbb{Q}[G]e_{i} \longrightarrow  M_{m_{i}}(\mathbb{Q}[N_{i}]\varepsilon_{i}) $$
\noindent given by
 $$ \alpha \stackrel{\theta_{i}}{\longmapsto}(\alpha_{rs})_{m_{i} \times m_{i}}.$$ \noindent As $\varepsilon_{i}^{t}$, $t \in T_{i}$, are mutually orthogonal idempotents and $\sum_{t \in T_{i}}\varepsilon_{i}^{t}=e_{i}$, it can be checked that $\theta_{i}$ is an isomorphism of $\mathbb{Q}$-algebras. This isomorphism in turn yields the group isomorphism $$\mathcal{Z}(\mathcal{U}( \mathbb{Q}[G]e_{i}))\stackrel{\theta_{i}}{\cong}\mathcal{Z}(\mathcal{U}(M_{m_{i}}(\mathbb{Q}[N_{i}]\varepsilon_{i})))$$

\noindent given by $$ \alpha \stackrel{\theta_{i}}{\longmapsto} \varepsilon_{i}\alpha \varepsilon_{i}I_{m_{i}}.$$ \noindent Let $\psi_{i}$ denote the canonical isomorphism from  $\mathcal{Z}(\mathcal{U}( \mathbb{Q}[N_{i}]\varepsilon_{i}))$ to $\mathcal{Z}(\mathcal{U}(M_{m_{i}} (\mathbb{Q}[N_{i}]\varepsilon_{i}))) $ given by $u\stackrel{\psi_{i}}{\mapsto}uI_{m_{i}}.$ Denote by $\phi_{i}$, the restriction of $\psi_{i}^{-1} o~ \theta_{i}$ to $\mathcal{Z}(\mathcal{U}( \mathbb{Z}[G]e_{i}))$. We assert that $\phi_{i}$ is an isomorphism from $\mathcal{Z}(\mathcal{U}( \mathbb{Z}[G]e_{i}))$ to $\mathcal{Z}(\mathcal{U}( \mathbb{Z}[N_{i}]\varepsilon_{i})).$ For this, we need to show that
\begin{equation}\label{eq3}
    \phi_{i}(\mathcal{Z}(\mathcal{U}( \mathbb{Z}[G]e_{i})))=\mathcal{Z}(\mathcal{U}( \mathbb{Z}[N_{i}]\varepsilon_{i})).
\end{equation}
\noindent Consider $\alpha = (\sum_{g \in G}\alpha_{g}g)e_{i} \in \mathcal{Z}(\mathcal{U}( \mathbb{Z}[G]e_{i}))$ with $\alpha_{g} \in \mathbb{Z}$. We have   $\varepsilon_{i}\alpha\varepsilon_{i} = \sum_{g \in N_{i}} \alpha_{g}g\varepsilon_{i} \in \mathbb{Z}[N_{i}]\varepsilon_{i}.$ Consequently, $\phi_{i}(\alpha)=\varepsilon_{i}\alpha\varepsilon_{i} \in \mathcal{Z}(\mathcal{U}( \mathbb{Z}[N_{i}]\varepsilon_{i}))$, as we already have $\phi_{i}(\alpha) \in \mathcal{Z}(\mathcal{U}( \mathbb{Q}[N_{i}]\varepsilon_{i})).$ On the other hand, to see that $\mathcal{Z}(\mathcal{U}( \mathbb{Z}[N_{i}]\varepsilon_{i}))$ is contained in $\phi_{i}(\mathcal{Z}(\mathcal{U}( \mathbb{Z}[G]e_{i}))) $, let $u\in \mathcal{Z}(\mathcal{U}( \mathbb{Z}[N_{i}]\varepsilon_{i})).$ Following the argument as in the proof of Theorem 3.5 of \cite{jesp}, it can be seen that $\sum_{t \in T_{i}}u^{t}$ belongs to $\mathcal{Z}(\mathcal{U}( \mathbb{Z}[G]e_{i}))$, as $\varepsilon^{t}$, $t \in T_{i}$, are mutually orthogonal idempotents. One checks that $\sum_{t \in T_{i}}u^{t}$ maps to $u$ under $\phi_{i}$ and hence Eq.(\ref{eq3}) follows.
  \noindent The isomorphism $\phi_{i}$ now provides the group isomorphism $$\Theta_{i} : \mathcal{A}(H_{i},K_{i}) \longrightarrow A(H_{i},K_{i})$$
 \noindent by setting $$ 1-e_{i} + \alpha  \stackrel{\Theta_{i}}{\mapsto}  1-\varepsilon_{i}+  \varepsilon_{i}\alpha   \varepsilon_{i},$$ \noindent where $\alpha \in \mathcal{Z}(\mathcal{U}( \mathbb{Z}[G]e_{i})).$ We further see that if $u =  1-\varepsilon_{i} + \gamma \in B(H_{i}, K_{i})$, with $ \gamma \in \mathcal{Z}(\mathcal{U}(\mathbb{Z}[N_{i}]\varepsilon_{i})),$ then $\Theta_{i}(\prod_{t\in T_{i}}u^{t}) = \Theta_{i}(1-e_{i} + \sum_{t\in T_{i}}\gamma^{t}) =  1-\varepsilon_{i}+  \varepsilon_{i} (\sum_{t\in T_{i}}\gamma^{t}) \varepsilon_{i}=  1-\varepsilon_{i} + \gamma= u$. This yields  $\Theta_{i}(\mathcal{B}(H_{i},K_{i})) = B(H_{i},K_{i})$ and consequently, Lemma \ref{l1} follows. $~\Box$\\

\begin{lemma}\label{l2} Let $p$ be a prime and let $n \geq 1$ be an integer. For a subgroup $A$ of $\mathcal{U}(\mathbb{Z}/p^{n}\mathbb{Z})(\cong Aut(\langle \zeta_{p^{n}}\rangle))$, let $\mathcal{U}(\mathbb{Z}[\zeta_{p^{n}}]^{A})$ denote the unit group of the fixed ring $\mathbb{Z}[\zeta_{p^{n}}]^{A}$. If $\langle A,-1 \rangle =\mathcal{U}(\mathbb{Z}/p^{n}\mathbb{Z}),$ then

$$\mathcal{U}(\mathbb{Z}[\zeta_{p^{n}}]^{A})=
\left\{
  \begin{array}{ll}
\langle \zeta_{4} \rangle, & \hbox{if $p=2$, $A=\langle r \rangle$, $r\equiv 1({\rm mod}~ 4)$, $n \geq 2$;} \\
    \langle \zeta_{6} \rangle, & \hbox{if $p=3$, $A=\langle r \rangle$, $r\equiv 1({\rm mod}~ 3)$;} \\
     \langle \zeta_{2} \rangle, & \hbox{otherwise.}
  \end{array}
\right.
$$

\end{lemma}

\noindent\textbf{Proof.} If $\langle A,-1 \rangle =\mathcal{U}(\mathbb{Z}/p^{n}\mathbb{Z}),$ then it follows from (\cite{jesp}, Lemma 3.2) that the rank of $\mathcal{U}(\mathbb{Z}[\zeta_{p^{n}}]^{A})=0$ and hence it consists only of the roots of unity. Recall that the group of roots of unity in $\mathbb{Q}(\zeta_{p^{n}})$ is $\langle \zeta_{p^{n}}, -1\rangle$. Thus, in case $\langle A,-1 \rangle =\mathcal{U}(\mathbb{Z}/p^{n}\mathbb{Z}),$ any $u \in \mathcal{U}(\mathbb{Z}[\zeta_{p^{n}}]^{A})$ is of the type $\zeta_{p^{n}}^{j}$ or $\pm \zeta_{p^{n}}^{j}$ depending on $p$ is even or odd with $0 \leq j < p^{n}$ satisfying
\begin{equation}\label{eqn 4}
    jx \equiv j ({\rm mod}~ p^{n})~ \forall~ x \in A.
\end{equation}
If $p$ is odd, we have $A=\langle r \rangle$ for some $r\in \mathcal{U}(\mathbb{Z}/p^{n}\mathbb{Z})$ and if $p=2$, then $A=\langle r \rangle$ or $\langle r \rangle \times \langle -1 \rangle$ for some $r\in \mathcal{U}(\mathbb{Z}/2^{n}\mathbb{Z})$.\\

\noindent\underline{Case(I): $A=\langle r \rangle$}\\

\noindent We must have $o_{p^{n}}(r)=\varphi(p^{n})$ or $\frac{\varphi(p^{n})}{2}$.\\

\noindent Suppose $o_{p^{n}}(r)=\varphi(p^{n})$. This case arises only if $p$ is an odd prime or $p^{n}=2$ or $4$. From Eq.(\ref{eqn 4}), $-j \equiv j({\rm mod}~ p^{n})$. This gives $j=0,~2^{n-1}$, if $p^{n}=2$ or $4$ and $j=0$, if $p$ is an odd prime.\\

\noindent If $o_{p^{n}}(r)=\frac{\varphi(p^{n})}{2}$, then Eq.(\ref{eqn 4}) holds if and only if
\begin{equation}\label{eqn 1}
    jr \equiv j ({\rm mod}~ p^{n}).
\end{equation}
Assuming $j \neq 0$, write $j=p^kj^{'}$ with $p \nmid j^{'}$ and $k \geq 0$. Eq.(\ref{eqn 1}) implies that $r^{p^{k}} \equiv 1({\rm mod}~p^{n})$ which gives that $\frac{\varphi(p^{n})}{2}$ divides $ ~p^{k}.$ This is possible only if $p=2$ or 3. Furthermore, for $p=2$, we have $n\geq 2$ and $j=2^{n-2},2.2^{n-2},3.2^{n-2}$; and for $p=3$, we have $j=3^{n-1},2.3^{n-1}.$\\

\noindent Consequently,
$$\mathcal{U}(\mathbb{Z}[\zeta_{p^{n}}]^{A}) = \left\{
                                     \begin{array}{ll}
                                       \{\zeta_{2^{n}}^{0},\zeta_{2^{n}}^{2^{n-2}},(\zeta_{2^{n}}^{2^{n-2}})^{2},(\zeta_{2^{n}}^{2^{n-2}})^{3}\}=\langle\zeta_{4}\rangle, & \hbox{if $r \equiv 1({\rm mod}~ 4)$,~$n \geq 2$;} \vspace{0.2cm}\\
                                       \{\pm\zeta_{3^{n}}^{0},\pm\zeta_{3^{n}}^{3^{n-1}},\pm(\zeta_{3^{n}}^{3^{n-1}})^{2}\}=\langle\zeta_{6}\rangle, & \hbox{if $r \equiv 1({\rm mod}~ 3)$;} \vspace{0.2cm}\\
                                       \{\pm 1\}=\langle\zeta_{2}\rangle, & \hbox{otherwise.}
                                     \end{array}
                                   \right.$$\\

\noindent\underline{Case(II): $A=\langle r \rangle \times \langle -1 \rangle $}\\

\noindent This case arises only if $p=2$, $n\geq 3$ and Eq.(\ref{eqn 4}) holds $\Leftrightarrow 2j\equiv 0 ({\rm mod}~2^{n}) \Leftrightarrow j=0,~2^{n-1}$. This yields the desired result. $~\Box$\\

\noindent{\bf Proof of Theorem \ref{t1}.}
(i) From (\cite{jesp}, Theorem 3.1), it follows immediately that rank of $\mathcal{Z}(\mathcal{U}(\mathbb{Z}[G]))=0$ if and only if $|I_{i}|=1$, $\forall ~i$, $1\leq i \leq m$. Further, (\cite{sehgal}, Corollary 7.3.3) implies that $\mathcal{Z}(\mathcal{U}(\mathbb{Z}[G]))=\pm \mathcal{Z}(G)$ in this case.\\

\noindent(ii) Since  $\mathcal{Z}(\mathcal{U}(\mathbb{Z}[G]))$ is a subgroup of $ \prod_{i=1}^{m}\mathcal{A}(H_{i},K_{i}),$ we have

 \begin{eqnarray}\label{e2}
                                            [\mathcal{Z}(\mathcal{U}(\mathbb{Z}[G])):\langle\mathcal{B}(G)\rangle] &\leq & [\prod_{i=1}^{m}\mathcal{A}(H_{i},K_{i}):\langle\mathcal{B}(G)\rangle] \nonumber\\
                                            &=& [\prod_{i=1}^{m}\mathcal{A}(H_{i},K_{i}):\langle\cup_{i=1}^{m}\mathcal{B}(H_{i},K_{i})\rangle] \nonumber\\
                                            &=& \prod_{i=1}^{m}[\mathcal{A}(H_{i},K_{i}):\langle\mathcal{B}(H_{i},K_{i})\rangle]\nonumber\\
                                            &=& \prod_{i=1}^{m}[A(H_{i},K_{i}):\langle B(H_{i},K_{i})\rangle].
                                         \end{eqnarray}

\noindent The last equality follows from Lemma \ref{l1}. We now show that
 \begin{equation}\label{eqn 6}
  [A(H_{1},K_{1}):\langle B(H_{1},K_{1})\rangle] \leq 2
\end{equation}

\noindent and for $2 \leq i \leq m$,

\begin{equation}\label{eqn 3}[A(H_{i},K_{i}):\langle B(H_{i},K_{i})\rangle] \leq \left\{
                                                        \begin{array}{ll}
                                                          o_{i}, & \hbox{if $|I_{i}|=1$;} \vspace{0.2cm}\\
                                                         2h_{p_{i}^{n_{i}}}^{+}l_{i}p_{i}^{n_{i}-1} \mathfrak{o}_{i} (l_{i}^{d_{i}-1}[N_{i} :H_{i}])^{|I_{i}|-1}, & \hbox{if $|I_{i}|\neq 1$.}
                                                        \end{array}
                                                      \right.
\end{equation}

\noindent Let $1 \leq i \leq m.$  We have that the center of $\mathbb{Q}[N_{i}]\varepsilon_{i}$ is equal to $(\mathbb{Q}[H_{i}]\varepsilon_{i})^{N_{i}/H_{i}}$, where $(\mathbb{Q}[H_{i}]\varepsilon_{i})^{N_{i}/H_{i}}$ denotes the fixed field under the action of $N_{i}/H_{i}$  on $\mathbb{Q}(\zeta_{p_{i}^{n_{i}}})\cong\mathbb{Q}[H_{i}]\varepsilon_{i}.$ Now, the center of $\mathbb{Q}(1-\varepsilon_{i})+\mathbb{Q}[N_{i}]\varepsilon_{i}$, which is equal to $\mathbb{Q}(1-\varepsilon_{i})+(\mathbb{Q}[H_{i}]\varepsilon_{i})^{N_{i}/H_{i}}$, is embedded inside the algebra $\mathbb{Q}[H_{i}]\hat{K_{i}} \oplus \mathbb{Q}(1-\hat{K_{i}})$, via the embedding
\begin{equation}\label{e4} r(1-\varepsilon_{i})+u\stackrel{\iota}{\mapsto} (r(1-\varepsilon_{i})+u)\hat{K_{i}}+r(1-\hat{K_{i}}),\end{equation} where $r \in \mathbb{Q}$ and $u \in (\mathbb{Q}[H_{i}]\varepsilon_{i})^{N_{i}/H_{i}}.$\\

\noindent As $H_{i}/K_{i}=\langle g_{i}K_{i}\rangle$, any element $x\hat{K_{i}}\in \mathbb{Q}[H_{i}]\hat{K_{i}}$ can be expressed as $x\hat{K_{i}}=\sum_{j=0}^{p_{i}^{n_{i}}-1}x_{j}g_{i}^{j}\hat{K_{i}}$ with $x_{j}\in \mathbb{Q}.$ Let $\pi$  denote the projection of $\mathbb{Q}[H_{i}]\hat{K_{i}} \oplus \mathbb{Q}(1-\hat{K_{i}})$ onto $\mathbb{Q}(\zeta_{p_{i}^{n_{i}}})$ under the isomorphism $\mathbb{Q}[H_{i}]\hat{K_{i}} \oplus \mathbb{Q}(1-\hat{K_{i}}) \stackrel{\tau}{\cong} \oplus_{k=0}^{n_{i}}\mathbb{Q}(\zeta_{p_{i}^{k}})\oplus\mathbb{Q}(1-\hat{K_{i}})$ given by \begin{equation}\label{e5}x\hat{K_{i}}+a(1-\hat{K_{i}})\stackrel{\tau}{\mapsto}(\sum_{j=0}^{p_{i}^{n_{i}}-1}x_{j}, \sum_{j=0}^{p_{i}^{n_{i}}-1}x_{j}\zeta_{p_{i}}^{j},...,\sum_{j=0}^{p_{i}^{n_{i}}-1}x_{j}\zeta_{p_{i}^{n_{i}}}^{j},
a(1-\hat{K_{i}})),\end{equation} where $x\hat{K_{i}}=\sum_{j=0}^{p_{i}^{n_{i}}-1}x_{j}g_{i}^{j}\hat{K_{i}}\in \mathbb{Q}[H_{i}]\hat{K_{i}}$ and $a\in \mathbb{Q}. $\\

 \noindent Now observe that $\pi \circ \tau \circ \iota$ is injective on $\mathcal{Z}(\mathbb{Q}(1-\varepsilon_{i})+\mathbb{Q}[N_{i}]\varepsilon_{i})$ and
   $ \pi \circ \tau \circ \iota(A(H_{i},K_{i}))$ is contained in $ \mathcal{U}(\mathbb{Z}[\zeta_{p_{i}^{n_{i}}}]^{N_{i}/H_{i}}).$ Hence,
\begin{eqnarray}\label{e10}
   [A(H_{i},K_{i}): \langle B(H_{i},K_{i})\rangle] &=& [\pi \circ \tau \circ \iota(A(H_{i},K_{i})):\pi \circ \tau \circ \iota(\langle B(H_{i},K_{i})\rangle)] \nonumber \\
   &\leq & [\mathcal{U}(\mathbb{Z}[\zeta_{p_{i}^{n_{i}}}]^{N_{i}/H_{i}}): \pi \circ \tau \circ \iota(\langle B(H_{i},K_{i})\rangle)]. \end{eqnarray}

 \noindent If $i=1$, i.e., $(H_{i},K_{i})=(G,G)$, then $[A(H_{i},K_{i}):\langle B(H_{i},K_{i})\rangle] \leq |\mathcal{U}(\mathbb{Z})|=2$. Thus Eq.(\ref{eqn 6}) holds.\\

\noindent If $2 \leq i \leq  m$ is such that $|I_{i}|=1$, then $B(H_{i},K_{i})$ is an empty set and therefore, $[\mathcal{U}(\mathbb{Z}[\zeta_{p_{i}^{n_{i}}}]^{N_{i}/H_{i}}): \pi \circ \tau \circ \iota(\langle B(H_{i},K_{i})\rangle)]=|\mathcal{U}(\mathbb{Z}[\zeta_{p_{i}^{n_{i}}}]^{N_{i}/H_{i}})|.$ We have using Lemma \ref{l2},
\begin{equation}\label{eqn 8}
    [A(H_{i},K_{i}):\langle B(H_{i},K_{i})\rangle] \leq o_{i},
\end{equation}
as $\langle N_{i}/H_{i}, -1 \rangle = \mathcal{U}(\mathbb{Z}/p_{i}^{n_{i}}\mathbb{Z})$ in this case.\\

\noindent We next assume that $|I_{i}|\not=1$.\\

 \noindent Set
   $$\begin{array}{cl}N_{(H_{i},K_{i})}= & \langle{\pi_{N_{i}/H_{i}}(\eta_{k}(\zeta_{p_{i}^{n_{i}}})^{o_{p_{i}^{n_{i}}}(k)p_{i}^{n_{i}-1}n_{H_{i},K_{i}}}})~|~ k \in I_{i}\setminus \{1\} \rangle,\\F_{(H_{i},K_{i})}= & \langle \eta_{k}(\zeta_{p_{i}^{n_{i}}})^{o_{p_{i}^{n_{i}}}(k)p_{i}^{n_{i}-1}n_{H_{i},K_{i}}}~|~ 1< k < \frac{p_{i}^{n_{i}}}{2}, (k,p_{i})=1 \rangle,
  \\O_{(H_{i},K_{i})}= & F_{(H_{i},K_{i})} \times \langle \zeta_{p_{i}^{n_{i}}}^{p_{i}^{n_{i}-1}}, -1 \rangle,
  \\ P_{(H_{i},K_{i})}= & \langle \eta_{k}(\zeta_{p_{i}^{n_{i}}})~|~ k \in \mathcal{U}(\mathbb{Z}/p_{i}^{n_{i}}\mathbb{Z}) \rangle
  = \langle  \eta_{k}(\zeta_{p_{i}^{n_{i}}})~|~ 1< k < \frac{p_{i}^{n_{i}}}{2}, (k,p_{i})=1 \rangle \times \langle \zeta_{p_{i}^{n_{i}}}, -1\rangle,
  \\ Q_{(H_{i},K_{i})}= & \mathcal{U}(\mathbb{Z}[\zeta_{p_{i}^{n_{i}}}]^{N_{i}/H_{i}})\cap O_{(H_{i},K_{i})},
  \end{array}
   $$
   \noindent where $\pi_{N_{i}/H_{i}}(u)= \prod_{\sigma \in N_{i}/H_{i}}\sigma(u),~{\rm for~} u \in \mathbb{Q}(\zeta_{{p_{i}}^{n_{i}}}).$\\

\noindent By (\cite{jesp}, Proposition 3.4), \begin{eqnarray}\label{e11}\pi \circ \tau \circ \iota(\langle B(H_{i},K_{i})\rangle) & = & N_{(H_{i},K_{i})}.  \end{eqnarray}
   \noindent Therefore, \begin{eqnarray}\label{e13}[\mathcal{U}(\mathbb{Z}[\zeta_{p_{i}^{n_{i}}}]^{N_{i}/H_{i}}): \pi \circ \tau \circ \iota(\langle B(H_{i},K_{i})\rangle)]& = & [\mathcal{U}(\mathbb{Z}[\zeta_{p_{i}^{n_{i}}}]^{N_{i}/H_{i}}):  N_{(H_{i},K_{i})}]\nonumber\\ & =& [\mathcal{U}(\mathbb{Z}[\zeta_{p_{i}^{n_{i}}}]^{N_{i}/H_{i}}):  Q_{(H_{i},K_{i})}][ Q_{(H_{i},K_{i})} :  N_{(H_{i},K_{i})}] \nonumber \\& \leq & [\mathcal{U}(\mathbb{Z}[\zeta_{p_{i}^{n_{i}}}]):  O_{(H_{i},K_{i})}][ Q_{(H_{i},K_{i})} :  N_{(H_{i},K_{i})}]. \end{eqnarray}
    Further, \begin{equation}\label{e13} [\mathcal{U}(\mathbb{Z}[\zeta_{p_{i}^{n_{i}}}]):  O_{(H_{i},K_{i})}]= [\mathcal{U}(\mathbb{Z}[\zeta_{p_{i}^{n_{i}}}]):  P_{(H_{i},K_{i})}][ P_{(H_{i},K_{i})} :  O_{(H_{i},K_{i})}]\end{equation}
Clearly,
\begin{equation}\label{e14} [P_{(H_{i},K_{i})} : O_{(H_{i},K_{i})}] =  p_{i}^{n_{i}-1}{\prod_{\stackrel{1< k < \frac{p_{i}^{n_{i}}}{2}}{(k,p_{i})=1}} o_{p_{i}^{n_{i}}}(k)p_{i}^{n_{i}-1}n_{H_{i},K_{i}}}=p_{i}^{n_{i}-1}\mathfrak{o}_{i}.\end{equation}
\noindent Also, by (\cite{kumm}, Theorem 8.2),
  \begin{equation}\label{e15}[\mathcal{U}(\mathbb{Z}[\zeta_{p_{i}^{n_{i}}}]): P_{(H_{i},K_{i})}] = h_{p_{i}^{n_{i}}}^{+}. \end{equation}

  \noindent Next, observe that $Q_{(H_{i},K_{i})}\cap F_{(H_{i},K_{i})}$ is a free abelian group, and by (\cite{jesp}, Lemma 3.2), it has rank at most $|I_{i}|-1$. Furthermore, any element of $ Q_{(H_{i},K_{i})}\cap F_{(H_{i},K_{i})}$ is of order at most $l_{i}^{d_{i}-1}|N_{i}/H_{i}|$ modulo $N_{(H_{i},K_{i})} \cap F_{(H_{i},K_{i})}$. To see this, let $u \in Q_{(H_{i},K_{i})}\cap F_{(H_{i},K_{i})}$ and write $u= \displaystyle\prod_{\stackrel{1< k < \frac{p_{i}^{n_{i}}}{2}}{(k,p_{i})=1}}(\eta_{k}(\zeta_{p_{i}^{n_{i}}})^{o_{p_{i}^{n_{i}}}(k)p_{i}^{n_{i}-1}n_{H_{i},K_{i}}})^{\alpha_{k}},~ \alpha_{k} \geq 0.$ Since $\pi_{N_{i}/H_{i}}(\eta_{r_{i}^{t}}(\zeta_{p_{i}^{n_{i}}}^{j}))=1$ and $\pi_{N_{i}/H_{i}}(\eta_{-j}(\zeta_{p_{i}^{n_{i}}}))=\pi_{N_{i}/H_{i}}(-\zeta_{p_{i}^{n_{i}}}^{-j})\pi_{N_{i}/H_{i}}(\eta_{j}(\zeta_{p_{i}^{n_{i}}}))$, for $i,j \geq 0$, it turns out that $u^{|N_{i}/H_{i}|l_{i}^{d_{i}-1}} = (\pi_{N_{i}/H_{i}}(u))^{l_{i}^{d_{i}-1}} \in N_{(H_{i},K_{i})} \cap F_{(H_{i},K_{i})}$.\\

   \noindent Consequently,  \begin{equation}\label{e16}
    [Q_{(H_{i},K_{i})}\cap F_{(H_{i},K_{i})} : N_{(H_{i},K_{i})}\cap F_{(H_{i},K_{i})}]\leq (l_{i}^{d_{i}-1}|N_{i}/H_{i}|)^{|I_{i}|-1}
\end{equation}
and therefore,
 \begin{eqnarray}\label{e17}
   [Q_{(H_{i},K_{i})}: N_{(H_{i},K_{i})}] &\leq& [Q_{(H_{i},K_{i})}: N_{(H_{i},K_{i})}\cap F_{(H_{i},K_{i})}] \nonumber\\
    &=& [Q_{(H_{i},K_{i})}: Q_{(H_{i},K_{i})}\cap F_{(H_{i},K_{i})}][Q_{(H_{i},K_{i})}\cap F_{(H_{i},K_{i})}: N_{(H_{i},K_{i})}\cap F_{(H_{i},K_{i})}] \nonumber \\
   &\leq& l_{i}(l_{i}^{d_{i}-1}|N_{i}/H_{i}|)^{|I_{i}|-1}.
 \end{eqnarray}

\noindent Finally, Eqs.(\ref{e10})-(\ref{e17}) yield Eq.(\ref{eqn 3}), which in view of Eq.(\ref{e2}) and Eq.(\ref{eqn 6})  complete the proof.     $~\Box$\\

It is known that if $G$ is an abelian group, then $\mathcal{U}(\mathbb{Z}[G])=\pm G$ if and only if $G$ is of exponent $1,~2, ~3,~ 4$ or $6$ (see \cite{units}, Theorem 2.7). We have the following:

\begin{cor}\label{c1}
Let $G$ be a strongly monomial group with a complete irredundant set $\{(H_{i},K_{i})~|~ 1\leq i\leq m\}$ of strong Shoda pairs such that  $\left[H_{i}:K_{i}\right]=1\,,2\,,3\,,4$ or $6$ for all $i$, $1\leq i\leq m$. Then, $\mathcal{Z}(\mathcal{U}(\mathbb{Z}[G]))=\pm \mathcal{Z}(G).$ \para In particular, if $G$ is a strongly monomial {\rm(}e.g. abelian by supersolvable{\rm)} group of exponent $1\,,2\,,3\,,4$ or $6$, then $\mathcal{Z}(\mathcal{U}(\mathbb{Z}[G]))=\pm \mathcal{Z}(G).$ However, the converse need not be true.

\end{cor}
\noindent{\bf Proof.} Here, $|I_{i}|=1$ for all $i$, $1\leq i\leq m$. Therefore, Theorem \ref{t1}(i) is applicable. The group $\mathcal{G}_{1}$ defined in section 4.1 for $p=3$ is an example of a strongly monomial group of exponent 9 satisfying $\mathcal{Z}(\mathcal{U}(\mathbb{Z}[G]))=\pm \mathcal{Z}(G).$ $\Box$\\

For abelian $p$-groups, Theorem \ref{t1} gives the following:

\begin{cor}\label{c2}Let $G$ be an abelian $p$-group, $p$ prime, and let $K_{i}$, $1\leq i\leq m$, be all the subgroups of $G$ with  cyclic quotient groups. Suppose $\left[G:K_{i}\right]= p^{n_{i}}$, for $1\leq i \leq m.$ Then, the rank of $\mathcal{U}(\mathbb{Z}[G])$ is non zero if and only if $p^{n_{i}}\geq 5$ for some $i$. In this case, the index of $\langle \mathcal{B}(G)\rangle$
 in $\mathcal{U}(\mathbb{Z}[G])$ is at most
$$2\displaystyle\prod_{i= 1}^{m}2^{n_{i}}h_{2^{n_{i}}}^{+}\big(\displaystyle\prod_{\stackrel{1< k< 2^{n_{i}-1}}{(k,2)=1}}2^{n_{i}}o_{2^{n_{i}}}(k)n_{G,K_{i}}\big), ~{\rm if}~ p=2.$$
and
$$\displaystyle\prod_{i=1}^{m}2p^{n_{i}}h_{p^{n_{i}}}^{+}\big(\displaystyle\prod_{\stackrel{1< k< \frac{p^{n_{i}}}{2}}{(k,p)=1}}2p^{n_{i}}o_{p^{n_{i}}}(k)n_{G,K_{i}}\big),~{\rm if}~ p\neq 2;$$

\noindent where an empty product equals $1$.

 \end{cor}

\section{Non Abelian groups of order $p^{n},~ n \leq 4$}
Let $G$ be a non abelian group of order $p^{n}$, $n\leq 4$. Observe that any such group, being metabelian, is normally monomial.

\subsection{Non Abelian groups of order $p^{3}$}

\noindent If $p=2$, then $G$ is either isomorphic to $D_{4}$, the dihedral group of order $8$ or is isomorphic to $Q_{8}$, the group of quaternions. Both these groups satisfy the hypothesis of (\cite{units}, Theorem 6.1). Therefore, we already know that the group of central units in the integral group ring of these groups consists of only the trivial units.\\

\noindent If $p$ is an odd prime, then $G$ is isomorphic to one of the following groups:
\begin{itemize}
  \item $\mathcal{G}_{1}= \langle a, b \,|\, a^{p^{2}}=b^{p}=1,~ ab= ba^{p+1}\rangle;$
  \item $\mathcal{G}_{2}= \langle a,b,c\, |\, a^{p}=b^{p}=c^{p}=1, ab=bac, ac=ca, bc=cb \rangle.$
\end{itemize}

\noindent In (\cite{BM}, Theorems 3 and 4), a complete and irredundant set of strong Shoda pairs of these groups has been found. Applying Theorem 3.1 of \cite{jesp}, we obtain that\\
$${\rm Rank ~of~} \mathcal{Z}(\mathcal{U}(\mathbb{Z}[\mathcal{G}_{i}]))=\frac{(p-3)(p+2)}{2},~i=1,2.$$\\
\noindent We now illustrate Theorem \ref{t1} in the particular cases, when $p=3$ or 5.\\

 \noindent $\underline{\textbf{p=3}}:$ In this case, the rank of $\mathcal{Z}(\mathcal{U}(\mathbb{Z}[\mathcal{G}_{i}]))=0$ and therefore, by Theorem \ref{t1}(i),  $\mathcal{Z}(\mathcal{U}(\mathbb{Z}[\mathcal{G}_{i}]))= \pm \mathcal{Z}(\mathcal{G}_{i}) ,~i=1,2.$\\

 \noindent  $\underline{\textbf{p=5}}:$ In this case, the rank of $\mathcal{Z}(\mathcal{U}(\mathbb{Z}[\mathcal{G}_{i}]))=7 ,~i=1,2.$ Using Remark \ref{r1}, we have computed the value of $n_{H,K}$ corresponding to each strong Shoda pair $(H,K)$ of the groups $\mathcal{G}_{1}$ and $\mathcal{G}_{2}$, which are tabulated below:

$$\begin{array}{ll}
    (H,K) & n_{H,K}\vspace{0.1cm}\\\hline
    &\\
(\mathcal{G}_{1},\mathcal{G}_{1}) & ~~~~1 \vspace{0.2cm}\\
    (\langle a^5,b \rangle,\langle b\rangle) & ~~~~5\vspace{0.2cm} \\
    (\mathcal{G}_{1},\langle a\rangle) & ~~~~5^{2} \vspace{0.2cm}\\
    (\mathcal{G}_{1},\langle a^{5},b\rangle) & ~~~~5^{2} \vspace{0.2cm}\\
    (\mathcal{G}_{1},\langle ab\rangle)& ~~~~5^{2} \vspace{0.2cm}\\
    (\mathcal{G}_{1}, \langle a^{2}b\rangle)& ~~~~5^{2} \vspace{0.2cm}\\
    (\mathcal{G}_{1},\langle a^{3}b\rangle)& ~~~~5^{2} \vspace{0.2cm}\\
    (\mathcal{G}_{1}, \langle a^{4}b\rangle)& ~~~~5^{2} \vspace{0.2cm}\\
\end{array}~~~~~~~~~~~~~~~~~~~~~~~~~~~~
  \begin{array}{ll}
    (H,K) & n_{H,K}\vspace{0.1cm}\\\hline
    &\\
(\mathcal{G}_{2}, \mathcal{G}_{2})&~~~~ 1 \vspace{0.2cm}\\
    (\langle a,c \rangle,\langle a\rangle) & ~~~~5\vspace{0.2cm} \\
    (\mathcal{G}_{2},\langle b,c\rangle) & ~~~~5^{2} \vspace{0.2cm}\\
    (\mathcal{G}_{2},\langle a,c\rangle) & ~~~~5^{2} \vspace{0.2cm}\\
    (\mathcal{G}_{2},\langle ab,c\rangle)& ~~~~5^{2} \vspace{0.2cm}\\
    (\mathcal{G}_{2},\langle a^{2}b,c\rangle)& ~~~~5^{2} \vspace{0.2cm}\\
    (\mathcal{G}_{2},\langle a^{3}b,c\rangle) & ~~~~5^{2} \vspace{0.2cm}\\
    (\mathcal{G}_{2},\langle a^{4}b,c\rangle) & ~~~~5^{2} \vspace{0.2cm}\\
  \end{array}$$
 \begin{small}
$$\begin{array}{cc}
 {\rm \textbf{Fig.1} ~Strong ~Shoda~ pairs ~of}~\mathcal{G}_{1} &~~~~~~~~~~~~~~~ {\rm \textbf{Fig.2}~ Strong ~Shoda~ pairs ~of}~\mathcal{G}_{2}\end{array}$$
\end{small}

\noindent  Theorem \ref{t1} and (\cite{kumm},~\S11.5) yield that $[\mathcal{Z}(\mathcal{U}(\mathbb{Z}[\mathcal{G}_{i}])):\langle \mathcal{B}(\mathcal{G}_{i})\rangle]\leq 2^{29}5^{27},~i=1,2$.

\subsection{Non Abelian groups of order $p^{4}$}
We first take the case, when $p=2$. Upto isomorphism, there are $9$ non isomorphic groups of order $2^{4}$ as listed in (\cite{burn}, \S 118). Except the following two groups:\begin{itemize}
                       \item $\mathcal{H}_{1}=\langle a, b : a^{8} =b^{2} =1, ba=a^{7}b  \rangle$;
                       \item $\mathcal{H}_{2}=\langle a, b : a^{8} =b^{4} =1, ba=a^{7}b, a^{4}=b^{2}  \rangle,$
                     \end{itemize}
the other non abelian groups of order $2^{4}$ again satisfy the hypothesis of (\cite{units}, Theorem 6.1). Hence, if $G$ is a non abelian group of order $2^{4}$ other than $\mathcal{H}_{1}$ and $\mathcal{H}_{2}$, then  $\mathcal{Z}(\mathcal{U}(\mathbb{Z}[G]))=\pm \mathcal{Z}(G).$

For the groups $\mathcal{H}_{1}$ and $\mathcal{H}_{2}$, we obtain using Theorem \ref{t3}, that $\{(\langle a \rangle,\langle 1\rangle), ~(\langle a \rangle,\langle a^{4}\rangle),$ $~(\mathcal{H}_{1},\langle a\rangle),(\mathcal{H}_{1},\langle a^{2},b\rangle),~ (\mathcal{H}_{1},\langle a^{2},ab\rangle) ,~ (\mathcal{H}_{1}, \mathcal{H}_{1})\}$  and   $\{(\langle a \rangle,\langle 1\rangle), (\langle a \rangle,\langle a^{4}\rangle), (\mathcal{H}_{2},\langle a\rangle),$ $  (\mathcal{H}_{2},\langle a^{2},b\rangle),(\mathcal{H}_{2},\langle a^{2},ab\rangle), (\mathcal{H}_{2}, \mathcal{H}_{2})\} $
  are complete irredundant sets of strong Shoda pairs of $\mathcal{H}_{1}$ and $\mathcal{H}_{2}$ respectively. Theorem 3.1 of \cite{jesp} now yields that $$ {\rm~ Rank ~of}~ \mathcal{Z}(\mathcal{U}(\mathbb{Z}[H_{i}]))= 1, ~i=1,2.$$

   \noindent Also, by Theorem \ref{t1} and (\cite{kumm}, \S 11.5), it follows that
 $$[\mathcal{Z}(\mathcal{U}(\mathbb{Z}[\mathcal{H}_{i}])):\langle \mathcal{B}(\mathcal{H}_{i})\rangle] \leq 2^{12},~i=1,2.$$
\noindent We next assume that $p$ is an odd prime.\\

 Upto isomorphism, the following are all the non abelian groups of order $p^{4}$ (see \cite{burn}, \S 117):

     \begin{enumerate}
       \item  $G_{1}= \langle a, b : a^{p^{3}}= b^{p}=1, ba=a^{1+p^{2}}b \rangle;$
       \item $G_{2}=\langle a, b, c: a^{p^{2}} =b^{p} =c^{p}=1, cb=a^{p}bc, ab=ba, ac=ca \rangle;$
       \item $G_{3}=\langle a, b : a^{p^{2}} =b^{p^2} =1, ba=a^{1+p}b  \rangle;$
       \item $G_{4}=\langle a, b, c: a^{p^{2}} =b^{p} =c^{p}=1, ca=a^{1+p}c, ba=ab, cb=bc \rangle;$
       \item $G_{5}=\langle a, b, c: a^{p^{2}} =b^{p} =c^{p}=1, ca=abc, ab=ba, bc=cb \rangle;$
       \item $G_{6}=\langle a, b, c: a^{p^{2}} =b^{p} =c^{p}=1, ba=a^{1+p}b, ca=abc, cb=bc \rangle;$
       \item $G_{7}=\begin{cases}\langle a, b, c: a^{p^{2}} =b^{p}=1, c^{p}=a^{p}, ab=ba^{1+p}, ac=cab^{-1}, cb=bc \rangle,~{\rm if}~ p=3, \\
                     \langle a, b, c: a^{p^{2}} =b^{p} =c^{p}=1, ba=a^{1+p}b, ca=a^{1+p}bc, cb=a^{p}bc \rangle, ~{\rm if}~ p>3 ;\end{cases}$
       \item $G_{8}=\begin{cases}\langle a, b, c: a^{p^{2}} =b^{p}=1, c^{p}=a^{-p}, ab=ba^{1+p}, ac=cab^{-1}, cb=bc \rangle, ~{\rm if}~ p=3, \\
\langle a, b, c: a^{p^{2}} =b^{p} =c^{p}=1, ba=a^{1+p}b, ca=a^{1+dp}bc, cb=a^{dp}bc \rangle,~ {\rm if}~p>3 \end{cases}$

\hspace{1.5cm} $ d\not\equiv 0,1({\rm mod}~ p);$
       \item $G_{9}=\langle a, b, c, d: a^{p} =b^{p} =c^{p}=d^{p}=1, dc=acd, bd=db, ad=da, bc=cb,\\~~~~~~~~~~~ac=ca,ab=ba \rangle;$
       \item $G_{10}=\begin{cases}\langle a, b, c: a^{p^{2}} =b^{p}=c^{p}=1, ab=ba, ac=cab, bc=ca^{-p}b \rangle,~ {\rm if}~ p=3,\\
                    \langle a, b, c,d: a^{p} =b^{p} =c^{p}=d^{p}=1, dc=bcd, db=abd,ad=da,\\~~~~~~~~~~~~~~bc=cb,ac=ca,ab=ba \rangle, ~{\rm if}~ p >3. \end{cases}$

     \end{enumerate}
\vspace{0.4cm}

    \begin{theorem}\label{t2}  The set $S(G_{i})$ of  strong Shoda pairs for each $G_{i}$, $1 \leq i \leq 10$, is as follows:\\
\begin{small} $$\begin{array}{cl}
  (i)~\mathcal{S}(G_{1})=&\hspace{-.4cm}\{ (\langle a \rangle,~\langle 1\rangle),~~ (G_{1},~\langle a\rangle),~~ (G_{1}, G_{1})\}~\cup \\
    &\hspace{-.4cm} \{(G_{1},~\langle a^{p^{2}},a^{pi}b\rangle),~~(G_{1},~\langle a^{p}, a^{i}b\rangle)~| ~0 \leq i \leq p-1 \};\vspace{0.45cm}\\
                  (ii)~ \mathcal{S}(G_{2})= &\hspace{-.4cm} \{(\langle a, b\rangle, ~\langle b \rangle),~~ (G_{2},~\langle a, b \rangle),~~(G_{2}, ~G_{2})\}~\cup \\
                    &\hspace{-.4cm} \{(G_{2},~\langle a,b^{i}c\rangle),~~ (G_{2},~\langle a^{i}b, a^{j}c\rangle)~| ~0 \leq i,j \leq p-1 \}
;\vspace{0.45cm}\\
(iii)~\mathcal{S}(G_{3})=  &\hspace{-.4cm} \{(G_{3},~\langle a,b^{p}\rangle),~~(G_{3},~\langle a\rangle),~~(G_{3}, ~G_{3})\}~\cup \\
   &\hspace{-.4cm} \{(\langle a,b^{p} \rangle,~\langle a^{pi}b^{p}\rangle),~~(G_{3},~\langle a^{p}, a^{i}b\rangle) ~|~0\leq i\leq p-1\} ~\cup\\
   &\hspace{-.4cm} \{(G_{3},~\langle a^{p}, a^{k}b^{p}\rangle)~|~ 1 \leq k\leq p-1\}
;\vspace{0.45cm}\\
(iv)~\mathcal{S}(G_{4})=  &\hspace{-.4cm} \{(G_{4},~\langle a, b \rangle),~~(G_{4}, ~G_{4})\} ~\cup\\
   &\hspace{-.4cm} \{(\langle a, b \rangle,~\langle a^{pi}b\rangle),~~(G_{4},~\langle a,b^{i}c\rangle),~~(G_{4},\langle a^{p},a^{i}b, a^{j}c\rangle)~|~0\leq i,j\leq p-1\}
;\vspace{0.45cm}\\
\end{array}$$
\end{small}

\begin{small} $$\begin{array}{cl}
(v)~\mathcal{S}(G_{5})= &\hspace{-.4cm} \{(\langle a, b \rangle,~\langle a \rangle),~~(G_{5},~\langle a^{p}, b,c \rangle),~~(G_{5},~\langle a,b\rangle),~~(G_{5}, ~G_{5})\}~\cup \\
   &\hspace{-.4cm} \{ (G_{5},~\langle b,a^{pi}c\rangle)~|~0 \leq i \leq p-1\}~\cup\\
   &\hspace{-.4cm} \{(\langle a, b \rangle,~\langle a^{p}b^{k}\rangle),~~(G_{5},\langle b,a^{k}c\rangle)~|~1 \leq k \leq p-1\}
;\vspace{0.45cm}\\
  (vi)~\mathcal{S}(G_{6})= &\hspace{-.4cm} \{(\langle a^{p}, b, c \rangle,~\langle a^{p}, c\rangle),~(G_{6},~\langle a, b \rangle),~(G_{6},~\langle a^{p},b,c\rangle, ),~(G_{6},~G_{6})\}~\cup \\
   &\hspace{-.4cm} \{(\langle a^{p}, b, c \rangle,~\langle b,a^{pi}c\rangle)~|~0 \leq i \leq p-1\}~\cup\\
   &\hspace{-.4cm}\{(G_{6},~\langle b, a^{k}c\rangle)~|~1 \leq k \leq p-1\}
;\vspace{0.45cm}\\
(vii)~\mathcal{S}(G_{7})= &\hspace{-.4cm} \{(\langle b,c \rangle, ~\langle  b\rangle),~~ (\langle b, c\rangle, ~\langle c \rangle),~~(G_{7},~\langle a,b\rangle),~~(G_{7}, ~G_{7})\}~\cup \\
  ~~~~(p=3 )& \hspace{-.4cm}\{(G_{7},~\langle  b,a^{i}c \rangle)~|~0 \leq i \leq p-1\}
 ;\vspace{0.45cm}\\
   (viii)~\mathcal{S}(G_{7})=&\hspace{-.4cm} \{(\langle b,ac\rangle, ~\langle  b \rangle),~~(\langle b, ac\rangle, ~\langle ac \rangle),~~(G_{7},~\langle a,b \rangle),~~(G_{7}, ~G_{7})\}~\cup \\
   \vspace{0.45cm}
~~~~(p>3)&\hspace{-.4cm} \{(G_{7},~\langle b, a^{i}c\rangle)~|~0 \leq i \leq p-1\}
;\vspace{0.45cm}\\
(ix)~\mathcal{S}(G_{8})= &\hspace{-.4cm} \{(\langle b,c\rangle, ~\langle  b\rangle),~~ (\langle b, c\rangle, ~\langle c \rangle),~~(G_{8},~\langle a,b\rangle),~~(G_{8}, ~G_{8})\}~\cup \\
  ~~~~(p=3 )& \hspace{-.4cm}\{(G_{8},~\langle  b,a^{i}c \rangle)~|~0 \leq i \leq p-1 \}
;\vspace{0.45cm}\\
   (x)~\mathcal{S}(G_{8})=&\hspace{-.4cm} \{ (\langle b, a^{d}c \rangle, ~\langle  b\rangle),~~ (\langle  b,a^{d}c\rangle, ~\langle  a^{d}c \rangle),~~ (G_{8},~\langle a,b \rangle),~(G_{8}, ~G_{8})\}~\cup \\
   \vspace{0.45cm}
~~~(p>3)&\hspace{-.4cm} \{(G_{8}, ~\langle  b,a^{i}c \rangle)~|~0 \leq i \leq p-1\}
;\vspace{0.45cm}\\
 (xi)~\mathcal{S}(G_{9})= & \hspace{-.4cm}\{(G_{9},~\langle a, b, d \rangle),~~(G_{9}, ~G_{9})\}~\cup \\
   &\hspace{-.4cm} \{(\langle a,b,d\rangle, ~\langle  d, a^{i}b \rangle),~~(G_{9},~\langle a, b, cd^{i}\rangle),~~(G_{9},~\langle a ,b^{i}c, b^{j}d\rangle)~|~0 \leq i,j \leq p-1\}
;\vspace{0.45cm}\\
 (xii)~\mathcal{S}(G_{10})= &\hspace{-.4cm} \{(\langle a, b\rangle, ~\langle a \rangle),~~(\langle  a,b\rangle, ~\langle  b \rangle),~~(G_{10},~\langle a,b\rangle),~~(G_{10}, ~G_{10})\}~\cup \\
  ~~~~(p=3 )&\hspace{-.4cm} \{(G_{10},~\langle  b,a^{i}c \rangle)~|~0 \leq i \leq p-1\};
 \vspace{0.45cm}\\
   (xiii)~\mathcal{S}(G_{10})=&\hspace{-.4cm} \{(\langle a,b,c\rangle, \langle  a,c \rangle),~~( G_{10},\langle a, b, d\rangle)\} ~\cup \\
   \vspace{0.45cm}
~~~(p>3)& \hspace{-.4cm} \{(\langle a,b,c \rangle,~\langle a^{i}c, b\rangle),~~(G_{10},\langle a,b,cd^{i} \rangle),~~(G_{10}, G_{10})~|~0 \leq i \leq p-1\}.
\end{array}$$
\end{small}
 \end{theorem}

\noindent{\bf Proof.} (i) Define $N_{0}:= \langle 1 \rangle$, $N_{1}:= \langle a^{p^{2}} \rangle$, $N_{2}:= \langle a^{p} \rangle$, $N_{3}:= \langle a \rangle$, $H_{i}:= \langle a^{p^{2}}, a^{pi}b \rangle$, $K_{j}:= \langle a^{p}, a^{j}b \rangle$ where $0 \leq i,j \leq p-1$. Observe that these subgroups are normal in $G_{1}$. Using Eq.(\ref{s1}), we have $\mathcal{S}_{N_{1}}= \mathcal{S}_{N_{2}}= \phi$, $\mathcal{S}_{N_{3}}=\{ (G_{1},N_{3} )\}$, $\mathcal{S}_{H_{i}}=\{ (G_{1},H_{i})\}$, $\mathcal{S}_{K_{j}}=\{ (G_{1},K_{j})\}$, $0 \leq i,j \leq p-1$. In order to find $\mathcal{S}_{N_{0}}$, we see that $\langle a \rangle$ is a maximal abelian subgroup of $G_{1}$. Further, the only subgroup $D$ of $\langle a \rangle$ satisfying $\operatorname{core}(D)=\langle 1\rangle$ is $D= \langle1 \rangle$. This gives $\mathcal{S}_{N_{0}}=\{ (\langle a \rangle,~\langle 1 \rangle)\}.$
Define $$\mathcal{N}_{1}= \{\langle 1 \rangle,~ \langle a^{p^{2}} \rangle,~\langle a^{p} \rangle,~ \langle a \rangle,~\langle a,b\rangle\} ~\cup~ \{\langle a^{p^{2}}, a^{pi}b \rangle,~ \langle a^{p}, a^{j}b \rangle~|~0 \leq i,j \leq p-1\}.$$ Observe that $\displaystyle \sum_{N\in \mathcal{N}_{1}}\sum_{D \in \mathcal{D}_{N}}[G:A_{N}]\varphi([A_{N}:D])=p^{4}$. Now, if $\mathcal{N}$ is the set of all normal subgroups of $G_{1}$, then \\

$ \begin{array}{cclr}
    p^{4} & =|G_{1}| &= \displaystyle\sum_{N\in \mathcal{N}}\sum_{D \in \mathcal{D}_{N}}[G:A_{N}]\varphi([A_{N}:D])&({\rm by ~Theorem ~ \ref{t3}}) \vspace{0.3cm}\\
&  &\geq\displaystyle \sum_{N\in \mathcal{N}_{1}}\sum_{D \in \mathcal{D}_{N}}[G:A_{N}]\varphi([A_{N}:D])&({\rm as}~ \mathcal{N}_{1}\subseteq \mathcal{N}) \vspace{0.3cm}\\
 &  &= p^{4}.\vspace{0.3cm} &
  \end{array}
$

\noindent This yields $\mathcal{S}_{N}= \phi$, if $N\not \in \mathcal{N}_{1}$ and consequently, by Theorem \ref{t3}, $\bigcup_{N\in \mathcal{N}_{1}}\mathcal{S}_{N}$ is a complete irredundant set of strong Shoda pairs of $G_{1}$.\\

\noindent (ii)-(xiii) For $2 \leq i \leq 10$, consider the following set $\mathcal{N}_{i}$ of normal subgroups of $G_{i}$:

\begin{small}
$$\begin{array}{cl}
\mathcal{N}_{2} ~=  &\{ \langle 1 \rangle,~\langle a^{p} \rangle, ~\langle a^{p},b \rangle,~\langle a,b \rangle,~\langle a,b,c \rangle \}~\cup\\&\{\langle a^{p}, b^{i}c \rangle,~\langle a, b^{i}c \rangle,~\langle ab^{i}c^{j} \rangle,~\langle a^{i}b, a^{j}c\rangle ~|~ 0\leq i,j \leq p-1 \};\vspace{0.3cm}\\
\mathcal{N}_{3} ~=  &\{ \langle 1 \rangle,~\langle a^{p} \rangle,~\langle a \rangle,~\langle a,b^{p} \rangle,~\langle a^{p},b^{p} \rangle, ~\langle a,b \rangle\}~\cup\\&\{~\langle a^{pi}b^{p} \rangle,~\langle a^{p}, a^{i}b \rangle ~|~ 0\leq i \leq p-1 \rangle\}~\cup ~\{ ~\langle a^{p},a^{k}b^{p} \rangle  ~|~ 1\leq k \leq p-1\};\vspace{0.3cm}\\
\mathcal{N}_{4} ~=  &\{ \langle 1 \rangle,~\langle a^{p} \rangle,~\langle a^{p},b \rangle,~\langle a,b \rangle,~\langle a,b,c \rangle \}~\cup\\&\{\langle a^{pi}b\rangle,~\langle a^{p}, b^{i}c \rangle,~\langle a, b^{i}c\rangle,~\langle ab^{i}c^{j} \rangle,~\langle  a^{p},a^{i}b, a^{j}c \rangle ~|~ 0\leq i,j \leq p-1 \rangle\};\vspace{0.3cm}\\
\mathcal{N}_{5} ~=  &\{ \langle 1 \rangle,~\langle b \rangle,~\langle a^{p},b \rangle,~\langle a^{p} \rangle,~\langle a^{p},b,c \rangle, ~\langle a, b \rangle,~\langle a,b,c \rangle\}~\cup\\&\{~\langle b,a^{pi}c \rangle ~|~ 0\leq i \leq p-1 \rangle\}~\cup ~\{ ~\langle a^{p}b^{k} \rangle,~\langle b,a^{k}c\rangle  ~|~ 1\leq k \leq p-1\};\vspace{0.3cm}\\
\mathcal{N}_{6} ~=  &\{\langle 1 \rangle,~\langle a^{p} \rangle, ~\langle a^{p},b \rangle,~\langle a,b \rangle,~\langle a^{p}, b,c \rangle,~\langle a,b,c \rangle\}~\cup ~\{ \langle b,a^{k}c \rangle  ~|~ 1\leq k \leq p-1\};\vspace{0.3cm}\\
\mathcal{N}_{7} ~=  &\{ \langle 1 \rangle,~\langle a^{p} \rangle,~\langle a^{p},b \rangle,~\langle a,b \rangle,~\langle a,b,c \rangle\}~\cup~ \{\langle b,a^{i}c \rangle, ~|~ 0\leq i \leq p-1 \};\vspace{0.3cm}\\
\mathcal{N}_{8} ~=  &\{ \langle 1 \rangle,~\langle a^{p} \rangle,~\langle a^{p},b \rangle,~\langle a,b \rangle,~\langle a,b,c \rangle\}~\cup ~\{\langle b,a^{i}c \rangle, ~|~ 0\leq i \leq p-1 \}\vspace{0.3cm};\\
\mathcal{N}_{9} ~=  &\{ \langle 1 \rangle,~\langle a \rangle,~\langle a,d \rangle,~\langle a,b,d \rangle , ~\langle a,b,c,d \rangle\}~\cup
\\&\{\langle a^{i}b \rangle,~\langle a,bc^{i}d^{j} \rangle,~\langle a,cd^{i} \rangle, ~\langle a,b,cd^{i} \rangle ~\langle a,b^{i}c,b^{j}d \rangle~|~ 0\leq i,j \leq p-1 \rangle\};\vspace{0.3cm}\\
\mathcal{N}_{10} ~= &\begin{cases}\{ \langle 1 \rangle,~\langle a^{3} \rangle, ~\langle a^{3},b \rangle,~\langle b,c \rangle,~\langle b,ac \rangle,~\langle b,a^{2}c \rangle,~\langle a, b \rangle, ~\langle a,b,c \rangle\}, ~{\mathrm if}~ p=3;\vspace{0.1cm}\\
                         \{
                        \langle 1 \rangle,~\langle a \rangle,~\langle a,b \rangle,~\langle a,b,d \rangle, ~\langle a,b,c,d \rangle \}~\cup~\{\langle a,b,cd^{i} \rangle~|~ 0\leq i \leq p-1 \}, ~{\mathrm if} ~p>3. \end{cases}
\end{array}$$
\end{small}

\vspace{0.5cm}

Now proceeding as in (i), we get the required complete and  irredundant set of strong Shoda pairs of $G_{i}$, $2 \leq i \leq 10$. $\Box$

\vspace{0.3cm}
\par For a particular odd prime $p$, the computation of $n_{H,K}$ corresponding to a strong Shoda pair $(H,K)\in \mathcal{S}(G_{i})$, $1\leq i \leq 10$, may be done using Remark \ref{r1}. An explicit bound on the index of $\langle\mathcal{B}(G_{i})\rangle$ in $\mathcal{Z}(\mathcal{U}(\mathbb{Z}[G_{i}]))$, $1\leq i \leq 10$, may thus be computed using Theorems \ref{t1} and \ref{t2}.

\begin{remark}\label{r2} It would be of interest to compute the integer $n_{H,K}$ corresponding to each strong Shoda pairs $(H,K)$ of the groups discussed in this section, explicitly in terms of $p$.\end{remark}

Finally, Theorem \ref{t2} along with (\cite{Oli}, Proposition 3.4) and (\cite{jesp}, Theorem 3.1) also yield the following:
\begin{cor} The Wedderburn decomposition of $\mathbb{Q}[G_{i}]$ and the rank of $\mathcal{Z}(\mathcal{U}(\mathbb{Z}[G_{i}]))$, $1\leq i \leq 10$, are as follows:
\begin{small}
\[\arraycolsep=1.4pt\def\arraystretch{2.2}
\begin{tabular}{|c|l|c|}
  \hline
  $G$ & $\mathbb{Q}[G]$ & \rm{Rank of} $\mathcal{Z}(\mathcal{U}(\mathbb{Z}[G]))$   \\ \hline\hline
   $G_{1}$ & $\mathbb{Q}\oplus \mathbb{Q}(\zeta_p)^{(1+p)}\oplus \mathbb{Q}(\zeta_{p^{2}})^{(p)}\oplus M_{p}(\mathbb{Q}(\zeta_{p^{2}}))$ & $\frac{(p+1)(p^{2}-5)}{2}$   \\  \hline
  $G_{2}$ & $\mathbb{Q}\oplus \mathbb{Q}(\zeta_p)^{(1+p+p^{2})}\oplus M_{p}(\mathbb{Q}(\zeta_{p^{2}}))$ & $\frac{p^{3}-p^{2}-3p-5}{2}$ \\\hline
  $G_{3}$ & $\mathbb{Q}\oplus \mathbb{Q}(\zeta_p)^{(1+p)}\oplus \mathbb{Q}(\zeta_{p^{2}})^{(p)}\oplus M_{p}(Q(\zeta_{p}))^{(p)}$ & $\frac{p^{3}+p^{2}-7p-3}{2}$  \\\hline
  $G_{4}$ & $\mathbb{Q}\oplus \mathbb{Q}(\zeta_p)^{(1+p+p^{2})}\oplus M_{p}(\mathbb{Q}(\zeta_p))^{(p)}$ & $\frac{(p-3)(p+1)^{2}}{2}$  \\\hline
  $G_{5}$ & $ \mathbb{Q}\oplus \mathbb{Q}(\zeta_p)^{(1+p)} \oplus \mathbb{Q}(\zeta_{p^{2}})^{(p)} \oplus M_{p}(\mathbb{Q}(\zeta_p))^{(p)}$
  & $\frac{p^{3}+p^{2}-7p-3}{2}$ \\\hline
  $G_{6}$ & $\mathbb{Q}\oplus \mathbb{Q}(\zeta_p)^{(1+p)}\oplus M_{p}(\mathbb{Q}(\zeta_p))^{(1+p)}$ & $(p-3)(p+1)$  \\\hline
  $G_{7}$ & $\mathbb{Q}\oplus \mathbb{Q}(\zeta_p)^{(1+p)}\oplus M_{p}(\mathbb{Q}(\zeta_p))\oplus M_{p}(\mathbb{Q}(\zeta_{p^{2}}))$ & $p^{2}-p-4$  \\\hline
  $G_{8}$ & $\mathbb{Q}\oplus \mathbb{Q}(\zeta_p)^{(1+p)}\oplus M_{p}(\mathbb{Q}(\zeta_p))\oplus M_{p}(\mathbb{Q}(\zeta_{p^{2}}))$ & $p^{2}-p-4$ \\\hline

  $G_{9}$ & $\mathbb{Q}\oplus \mathbb{Q}(\zeta_p)^{(1+p+p^{2})}\oplus M_{p}(\mathbb{Q}(\zeta_p))^{(p)}$ & $\frac{(p-3)(p+1)^{2}}{2}$ \\\hline
  $G_{10}$& \multirow{2}*{$\mathbb{Q}\oplus \mathbb{Q}(\zeta_3)^{(4)}\oplus M_{3}(\mathbb{Q}(\zeta_3))\oplus M_{3}(\mathbb{Q}(\zeta_9))$} & \multirow{2}*{$2$}\vspace{-0.5cm} \\
  $(p=3)$  & & \\ \hline
  $G_{10}$& \multirow{2}*{$\mathbb{Q}\oplus \mathbb{Q}(\zeta_p)^{(1+p)}\oplus M_{p}(\mathbb{Q}(\zeta_p))^{(1+p)}$}& \multirow{2}*{$(p-3)(p+1)$}\vspace{-0.5cm}\\
  $(p>3)$ &&\\ \hline
\end{tabular}\]\end{small}
\end{cor}

\vspace{0.7cm}

\noindent\textbf{Acknowledgements}\\

\noindent The authors are grateful to the anonymous referees for their valuable comments and suggestions which has helped to write the paper in the present form.
\bibliographystyle{amsplain}
\bibliography{ReferencesBM2Rev}
\end{document}